\documentclass[12pt,a4paper]{amsart}
\usepackage{amscd}
\usepackage{verbatim}
\usepackage[dvips]{graphicx}
\theoremstyle{plain}
  \newtheorem{thm}{Theorem}[section]

  \newtheorem{conj}[thm]{Conjecture}
  
  \newtheorem{obs}[thm]{Observation}
\theoremstyle{definition}

\theoremstyle{remark}
  
  \newtheorem*{ack}{Acknowledgments}
\newcommand{\Z}{\mathbb{Z}}
\newcommand{\R}{\mathbb{R}}
\newcommand{\C}{\mathbb{C}}
\newcommand{\cs}{\operatorname{cs}}
\newcommand{\CS}{\operatorname{CS}}
\newcommand{\Tr}{\operatorname{Tr}}
\newcommand{\vol}{\operatorname{vol}}

\newcommand{\li}{{\operatorname{Li}}_2}
\newcommand{\I}{\sqrt{-1}}
\allowdisplaybreaks
\begin{document}
\author
{Hitoshi Murakami}
\address
[Hitoshi Murakami]
{Department of Mathematics, Tokyo Institute of Technology,
 Oh-okayama, Meguro, Tokyo 152-8551, Japan}
\email
{starshea@tky3.3web.ne.jp}
\author
{Jun Murakami}
\address
[Jun Murakami]
{Department of Mathematical Sciences, 
School of Science and Engineering
Waseda University
3-4-1, Ohkubo
Shinjuku-ku, Tokyo
169-8555 JAPAN
}
\email
{murakami@mn.waseda.ac.jp}
\author
{Miyuki Okamoto}
\address
[Miyuki Okamoto]
{University of the Sacred Heart,
 Hiroo, Shibuya, Tokyo 150-8938, Japan}
\email
{miyuki3@hh.iij4u.or.jp}
\author
{Toshie Takata}
\address
[Toshie Takata]
{Department of Mathematics, Faculty of Science, Niigata University,
Niigata 950-2181, Japan}
\email
{takata@math.sc.niigata-u.ac.jp}
\author
{Yoshiyuki Yokota}
\address
[Yoshiyuki Yokota]
{Department of Mathematics, Tokyo Metropolitan University,
 Tokyo 192-0397, Japan}
\email
{jojo@math.metro-u.ac.jp}
\title
{Kashaev's conjecture and the Chern--Simons invariants of knots and links}
\begin{abstract}
R.M.~Kashaev conjectured that the asymptotic behavior of his link invariant,
which equals the colored Jones polynomial evaluated at a root of unity,
determines the hyperbolic volume of any hyperbolic link complement.
We observe numerically that for knots $6_3$, $8_9$ and $8_{20}$ and for
the Whitehead link, the colored Jones polynomials  are related to 
the hyperbolic volumes and the Chern--Simons invariants and propose a
complexification of Kashaev's conjecture.
\end{abstract}
\thanks
{This research is partially supported by Grand-in-Aid for Scientific
Research, The Ministry of Education, Science, Sports and Culture.}
\maketitle

%%%%%%%%%%%%%%%%%%%%%%%%%%%%%%%%%%%%%%%%%%%%%%%%%%%%%%%%%%%%%%%%%%%%%
\section{introduction}
%%%%%%%%%%%%%%%%%%%%%%%%%%%%%%%%%%%%%%%%%%%%%%%%%%%%%%%%%%%%%%%%%%%%%

In \cite{Kashaev:MODPLA95}, R.M.~Kashaev defined a link invariant 
associated with the quantum dilogarithm, depending on a positive 
integer $N$, which is denoted by $\langle L \rangle_N$ for a link $L$.
Moreover, in \cite{Kashaev:LETMP97}, he conjectured that for any 
hyperbolic link $L$, the asymptotics at 
$N\to \infty$ of $\left|\langle L \rangle_N \right|$ 
gives its volume, that is 
\begin{equation*}
  \vol(L)
  =
  2 \pi \lim_{N\to \infty} \frac{\log\left|\langle L \rangle_N\right|}{N}
\end{equation*}
with $\vol(L)$ the hyperbolic volume of the complement of $L$.
He showed that this conjecture is true for three doubled knots $4_1$, 
$5_2$, and $6_1$.
Unfortunately his proof is not mathematically rigorous.

Afterwards, in \cite{Murakami/Murakami:volume}, the first two authors 
proved that for any link $L$, Kashaev's invariant $\langle L \rangle_N$ 
is equal to the colored Jones polynomial evaluated at
$\exp\left(2\pi \sqrt{-1}/N\right)$, which is written by $J_N(L)$, and 
extended Kashaev's conjecture as follows.
\begin{conj}
[Volume Conjecture]
\begin{equation*}
  \Vert L \Vert
  =
  \frac {2\pi} {v_3} \lim_{N\to \infty} \frac {\log |J_N(L)|} {N}, 
  \label{eq:vol}
\end{equation*}
\noindent
where $\Vert L \Vert$ is the simplicial volume of the complement of $L$ 
and $v_3$ is the volume of the ideal regular tetrahedron.
\end{conj}

\noindent 
Note that the hyperbolic volume $\vol(L)$ of a hyperbolic link $L$ is equal
to $\Vert L \Vert$ multiplied by $v_3$.
This conjecture is not true for links in general, as $J_N(L)$ vanishes for 
a split link $L$.
Note also that it is shown by Kashaev and O. Tirkkonen in 
\cite{Kashaev/Tirkkonen:1999} that the volume conjecture holds for torus 
knots.
See \cite{D.Thurston:Grenoble} and \cite{Yokota:Murasugi70,Yokota:volume} 
for discussions about Kashaev's conjecture for hyperbolic knots from the 
viewpoint of tetrahedron decomposition.

In this paper, following Kashaev's way to analyze the asymptotic behavior
of the invariant, we observe numerically, by using MAPLE V (a product of 
Waterloo Maple Inc.) and SnapPea \cite{Weeks:SnapPea}, that for the 
hyperbolic knots $6_3$, $8_9$, $8_{20}$, and for the Whitehead link, the 
colored Jones polynomials  are related to the hyperbolic volumes and the 
Chern--Simons invariants.
Note that the knots $6_3$ and $8_9$ are not doubles of the unknot.

We also discuss a relation between the asymptotic behavior of $J_N(L)$ and 
the Chern--Simons invariant of the complement of the above-mentioned links 
$L$, and propose the following conjecture.
%
\begin{comment}

\begin{conj}[Complexification of Kashaev's conjecture]
Let $L$ be a hyperbolic link. 
Then the following formula holds.
  \begin{equation*}
    2 \pi \lim_{N\to\infty} \frac{\log{J_N(L)}}{N}
    =
    \vol(L)+\sqrt{-1}\CS(L) \mod \pi^2 \sqrt{-1},
  \end{equation*}
where $\CS(L)$ is the Chern--Simons invariant of $L$ 
\cite{Chern/Simons:ANNMA274,Meyerhoff:density}.
Note that the complement of $L$ is a hyperbolic manifold with cusps.
\end{conj}
\end{comment}
%
%
\begin{conj}[Complexification of Kashaev's conjecture]
Let $L$ be a hyperbolic link. 
Then the following formula holds.
  \begin{equation*}
    J_N(L) \sim \exp \frac {N }{2 \pi} (\vol(L)+\sqrt{-1}\CS(L))
    \quad (N\to \infty)   
  \end{equation*}
where $\CS(L)$ is the Chern--Simons invariant of $L$ 
\cite{Chern/Simons:ANNMA274,Meyerhoff:density}.
Note that the complement of $L$ is a hyperbolic manifold with cusps.
\end{conj}

 The statement of this conjecture  will be given more properly 
in the last section.
\begin{ack}
We thank the participants in the meeting
``Volume conjecture", October 1999 and those in the workshop
``Recent  Progress toward the Volume Conjecture", March 2000,
both of which were held at the International Institute for Advanced Study.
The latter was financially supported by the Research Institute for 
Mathematical Sciences, Kyoto University.
We are grateful to both of the institutes.

H.M., J.M. and M.O. express their gratitude to Graduate School of 
Mathematics, Kyushu University, where the essential part of this work was 
carried out in December 1999.

Thanks are also due to S.~Kojima for his suggestion of the Chern--Simons 
invariant,  to K.~Mimachi for valuable discussions, and 
to K.~Hikami for informing Pari-Gp \cite{Cohen:parigp} and fitting.
\end{ack}

%%%%%%%%%%%%%%%%%%%%%%%%%%%%%%%%%%%%%%%%%%%%%%%%%%%%%%%%%%%%%%%%%%%%%
\section{Preliminaries}
%%%%%%%%%%%%%%%%%%%%%%%%%%%%%%%%%%%%%%%%%%%%%%%%%%%%%%%%%%%%%%%%%%%%%

First we will briefly review the colored Jones polynomials of links 
following \cite{Kirby/Melvin:INVEM91}.
It is obtained from the quantum group $U_q(sl(2,\C))$  and its 
$N$-dimensional irreducible representation.

Let $L$ be an oriented link.
We consider a $(1,1)$-tangle presentation of $L$, obtained by cutting a 
component of the link.
We assume that all crossing and local extreme points are as in Figure~1.
We can calculate the $N$-colored Jones polynomial $J_{L}(N)$ evaluated at 
the $N$-th root of unity for $L$ in the following way.
We start with a labeling of the edges of the $(1,1)$-tangle presentation
with labels $\{0,1,\dots,N-1\}$.
Here we label the two edges containing the end points of the tangle by $0$.
Following the labeling, we associate a positive (respectively negative) 
crossing with the element  $R_{kl}^{ij}$ (respectively 
${\bar R}_{kl}^{ij}$), a maximal point $\cap$ labeled by $i$ with the 
element $-s^{-2i-1}$, and a minimal point $\cup$ labeled by $i$ with the 
element $-s^{2i+1}$ with $s=\exp\left(\frac {\pi \sqrt{-1}} N\right)$ as 
in Figure~1.

\vspace{-5mm}

\hspace*{13mm}
%%%%%%%%%%%%%%%%%%%%%
\includegraphics[scale=0.8]{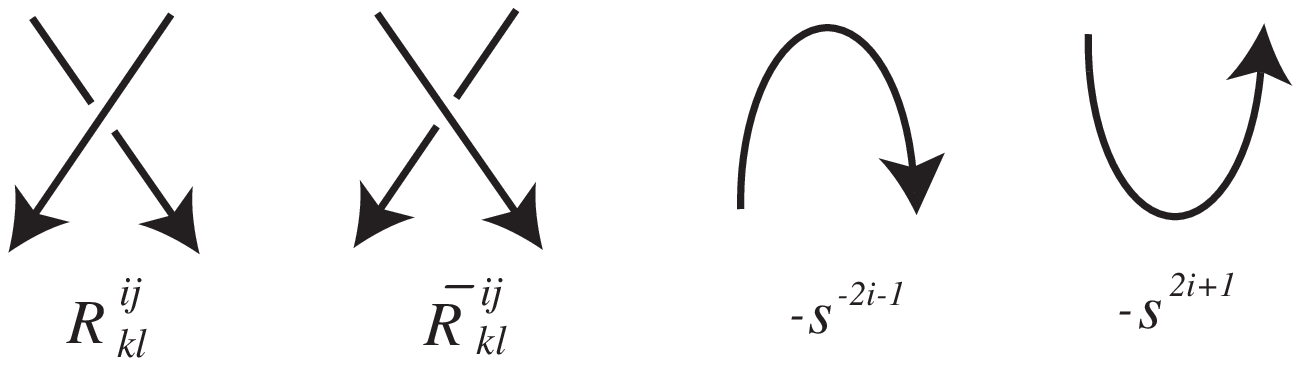}
%%%%%%%%%%%%%%%%%%%%%

\vspace{2mm}

\begin{center}
  Figure 1
\end{center}

\vspace{2mm}

\noindent
Here $R_{kl}^{ij}$ and ${{\bar R}}_{kl}^{ij}$ are given by 
\begin{align*}
  R_{kl}^{ij}
    & = \sum_{n=0}^{\min(N-1-i,j)}
        \delta_{l,i+n} \delta_{k,j-n} 
        \frac{(i+n)!(N-1+n-j)!}{(i)!(N-1-j)!(n)!}
        s^{2(i-\frac {N-1} 2)(j-\frac {N-1} 2)-n(i-j)-\frac {n(n+1)} 2},  \\
  {{\bar R}}_{kl}^{ij}
    & = \sum_{n=0}^{\min(N-1-j,i)}
        \delta_{l,i-n} \delta_{k,j+n} 
        \frac{(j+n)!(N-1+n-i)!}{(j)!(N-1-i)!(n)!}(-1)^n  \\
    & \hskip 3truecm
        \times
        s^{-2(i-\frac {N-1} 2)(j-\frac {N-1} 2)-n(i-j)+\frac {n(n+1)} 2}
\end{align*}
with $(n)!=(s-s^{-1})(s^2-s^{-2})\cdots (s^n-s^{-n})$.

After multiplying all elements associated to the critical points, we sum 
up over all labelings.
Here we ignore framings of links. 

Let us calculate the colored Jones polynomial of the Whitehead link as an 
example. 
We can label each edge in the following way, noting Kronecker's deltas in 
$R_{kl}^{ij}$ and ${\bar R}_{kl}^{ij}$.

\vspace{-6mm}

%%%%%%%%%%%%%%%%%%%%%
\hspace*{25mm}
\includegraphics[scale=1]{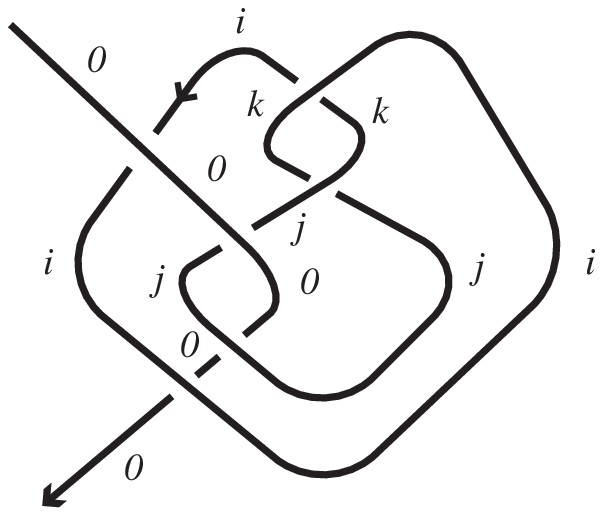}
%%%%%%%%%%%%%%%%%%%%%

\vspace{2mm}

\begin{center}
  Figure 2
\end{center}

\vspace{2mm}

We have to rotate a crossing where edges go up. In that case we use $\cup$ 
and/or $\cap$ to calculate the invariant.

Then we calculate the formula 
\begin{equation}\label{eq:Wh}
 J_N(L) = \sum_{\substack{0\le i,j,k\le N-1 \\ i,j\ge k}}
          \frac {(q)_i (q)_j \{(q)_{N-1-k}\}^2}
                {\{(q)_k\}^2 (q)_{N-1-i} (q)_{N-1-j} (q)_{i-k} (q)_{j-k}}
                 q^{-k(i+j+1)},
\end{equation}
where $q=s^2=\exp \left(\frac {2 \pi \sqrt{-1}} N\right)$.
Here $(x)_k=(1-x)(1-x^2)\cdots (1-x^k)$.

Next the Chern--Simons invariant of a link is defined as follows.
Let $\mathcal{A}$ be the set of all $SO(3)$-connections of the trivial
$SO(3)$-bundle of a closed three-manifold $M$ and 
$\cs \colon \mathcal{A}\to\R$ the Chern--Simons functional defined by
\begin{equation*}
  \cs(A)
  =
  \frac{1}{8\pi^2}\Tr
  \left(A \wedge d\,A+\frac{2}{3}A \wedge A \wedge A\right).
\end{equation*}
The Chern--Simons invariant of the connection $A$ is then defined to be
the integral
\begin{equation*}
  \cs_{M}(A) = \int_{s(M)}\cs(A)\in\R/\Z,
\end{equation*}
where the integral is over a section $s$ of the $SO(3)$-bundle (i.e., an 
orthonormal frame field on $M$) \cite{Chern/Simons:ANNMA274}.
If $M$ is hyperbolic we define $\cs(M)$ to be the Chern--Simons invariant
of the connection defined by the hyperbolic metric.

The definition of the Chern--Simons invariant for hyperbolic 
three-manifolds with cusps is due to R.~Meyerhoff \cite{Meyerhoff:density}.
It is defined modulo $1/2$ by using a special singular frame field which 
is linear near the cusps.
See \cite{Meyerhoff:density} for details.
See also \cite{Coulson/Goodman/Hodgson/Neumann:snap} how it is computed
by SnapPea \cite{Weeks:SnapPea}.
Throughout this paper we use another normalization 
$\CS(M)=-2\pi^2\cs(M)$ %$\CS(M)=2\pi^2\cs(M)$ 
so that $\vol(M)+\sqrt{-1}\CS(M)$ is a natural complexification of the 
hyperbolic volume $\vol(M)$ 
(see \cite{Neumann/Zagier:TOPOL85,Yoshida:INVEM85}).

%%%%%%%%%%%%%%%%%%%%%%%%%%%%%%%%%%%%%%%%%%%%%%%%%%%%%%%%%%%%%%%%%%%%%
\section{knot $6_3$}
%%%%%%%%%%%%%%%%%%%%%%%%%%%%%%%%%%%%%%%%%%%%%%%%%%%%%%%%%%%%%%%%%%%%%

Let us calculate the colored Jones polynomial of the knot $6_3$ using
the labeling as in Figure~3. 

\vspace{2mm}

%%%%%%%%%%%%%%%%%%%%%
\hspace*{-13mm}
\includegraphics[scale=1]{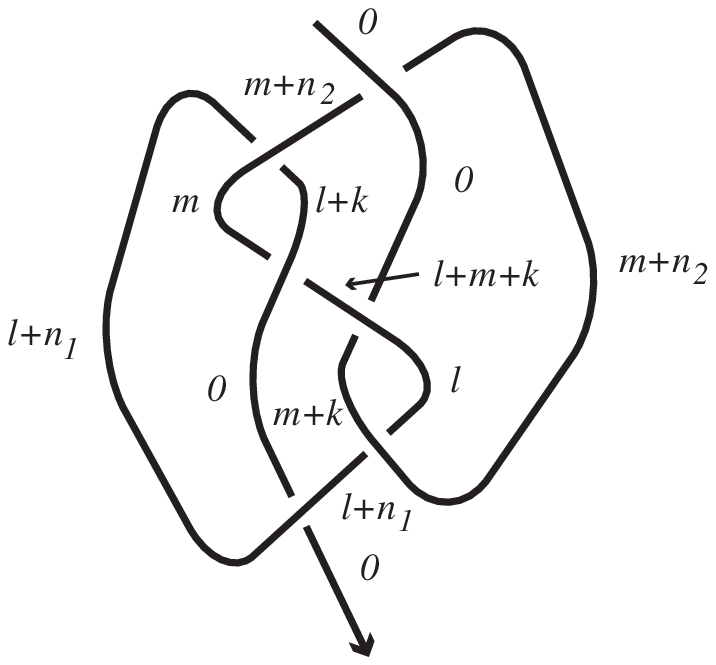}
%%%%%%%%%%%%%%%%%%%%%

\vspace{2mm}

\begin{center}
 Figure 3
\end{center}

\vspace{2mm}

\noindent
Putting $k=n_1+n_2$ and using the formula in \cite{Murakami/Murakami:volume}
\begin{equation}
  \sum_{i=0}^{N-1} (-1)^i s^{\beta i} 
  \begin{bmatrix} \alpha \\ i \end{bmatrix}
  =
  \prod_{j=1}^{\alpha} (1-s^{\beta+\alpha+1-2j}) \label{eq:knot}
\end{equation}
with $\alpha=k$, $i=n_1$, and $\beta=-k-1-2N$, we calculate 
\begin{align*}
  J_N (6_3) 
  =
  & \sum_{\scriptstyle  0 \le k,l,m 
          \atop
          \scriptstyle k+l+m \le N-1} 
    (-1)^{k+l}
    s^{\frac{(l+k)(l+k+1)}{2}
       - \frac{(m+k)(m+k+1)}{2}
       + \frac{k(k+1)}{2}
       + 2(m-l)(k+1)
       + N(m-l+k)}\\
  & \times
    \frac{(N-1-l)!(N-1-m)!(l+m+k)! (N-1)!(1-s^{-2N-2})
          \cdots (1-s^{-2N-2k})}
         {(N-1-l-m-k)! (N-1-l-k)! (N-1-m-k)! (l)!(m)!(k)!}.
\end{align*}

The colored Jones polynomial of the knot $6_3$ is given by 
\begin{equation}
  J_N (6_3)
  =
  \sum_{\substack{k,l,m\ge 0 \\ k+l+m \le N-1}}
  \left| \frac{(q)_{k+l+m}}{(q)_l (q)_m} \right|^2
  (q)_{k+l} (\bar q)_{m+k} \, q^{(m-l)(k+1)}.\label{eq:knot63}
\end{equation}

%%%%%%%%%%%%%%%%%%%%%%%%%%%%%%%
% Kashaev's way
%%%%%%%%%%%%%%%%%%%%%%%%%%%%%%%

We review of the technique in \cite{Kashaev:LETMP97}.
For a complex number $p$ and a positive real number $\gamma$ with
$|\operatorname{Re} p |<\pi + \gamma$, we define 
\begin{equation*}
  S_{\gamma}(p)
  =
  \exp \frac 1 4 \int_{-\infty}^\infty
  \frac{e^{px}}{\sinh (\pi x) \sinh (\gamma x)}
  \frac{dx}{x}.
\end{equation*}
Here $\operatorname{Re}$ denotes the real part.
This function has two properties:

\vspace{2mm}

(a)
\hspace{2mm}
$\displaystyle{
  \left(1+\exp(\sqrt{-1}p)\right) S_{\gamma}(p+\gamma)
  = S_{\gamma}(p-\gamma);
}$

\vspace{2mm}

(b)
\hspace{2mm}
$\displaystyle{
  S_{\gamma}(p) \sim
  \exp\left(
        \frac{1}{2\gamma \sqrt{-1}}\li \left(-\exp(\sqrt{-1}p) \right)
      \right)
  \hspace{10mm}
  (\gamma \to 0),
}$

\vspace{2mm}

\noindent
where 
\begin{equation*}
  \li(z) = -\int_0^z \frac{\log(1-u)}{u} du.
\end{equation*}
We put 
\begin{equation*}
  f_{\gamma}(p) = \frac{S_\gamma(\gamma - \pi)}{S_\gamma(p)},
  \quad
  \bar{f}_{\gamma}(p) = \frac{S_\gamma(-p)}{S_\gamma(\pi -\gamma)},
\end{equation*}
so that
\begin{equation*}
  (q)_{k} = f_{\gamma}(-\pi+(2k+1)\gamma), 
  \quad
  (\bar q)_k = \bar{f}_{\gamma}(-\pi+(2k+1)\gamma).
\end{equation*}

%%%%%%%%%%%%%%%%%%%%%%%%%%%%%%%%%%%%%%%%%%%%%%%%%%%%%%%%%%%%%%%%%%%

Following Kashaev's way, we rewrite the formula~\eqref{eq:knot63} as a
multiple integral with appropriately chosen contours.
(Note that there is considerable doubt as to the contours.)
By using the property (b), it can be asymptotically approximated by 
\begin{equation*}
  \int \! \! \! \int \! \! \! \int
  \exp \frac{\sqrt{-1}}{2\gamma}
  V_{6_3}(z,u,v) \, dz \, du \, dv
\end{equation*}
with $\gamma=\pi/N$.
Here $z$, $u$, and $v$ correspond to $q^k$, $q^m$, and $q^l$ respectively, 
and
\begin{align*}
  V_{6_3}(z,u,v)
  & = \li (zuv) - \li \left( \frac{1}{zuv} \right)
      + \li (zv) - \li \left( \frac{1}{zu} \right)  \\[2mm]
  & \hspace*{4mm}
      - \li (u) + \li \left( \frac{1}{u} \right)
      - \li (v) + \li \left( \frac{1}{v} \right)
      - \log z \log \frac{u}{v}.
\end{align*}
Then there exists a stationary point 
$$(z_0, u_0, v_0)=(0.204323-0.978904\sqrt{-1},
                  1.60838+0.558752\sqrt{-1},
                  0.554788+0.192734\sqrt{-1})$$
 of $V_{6_3}$ with
$$
  \operatorname{Im} V_{6_3} (z_0, u_0, v_0) < 0,
  \quad
  \arg z_0 + \arg u_0 + \arg v_0 \le 2 \pi,
$$
and we have 
\begin{align*}
  & - \operatorname{Im} V_{6_3}(z_0,u_0,v_0)
    = 5.693021 \ldots ,  \\[0.5mm]
  & \operatorname{Re} V_{6_3}(z_0,u_0,v_0) = 0.
\end{align*}
From values of $\vol(6_3)$ and $\CS(6_3)$ given by SnapPea, we see that
the equation
$$
  \exp \frac{\sqrt{-1}}{2\gamma} V_{6_3}(z_0,u_0,v_0)
  = \exp \frac{\vol(6_3) + \sqrt{-1} \CS(6_3)}{2\gamma}
$$
holds up to digits shown above.

%%%%%%%%%%%%%%%%%%%%%%%%%%%%%%%%%%%%%%%%%%%%%%%%%%%%%%%%%%%%%%%%%%%%%
\section{Knot $8_9$}
%%%%%%%%%%%%%%%%%%%%%%%%%%%%%%%%%%%%%%%%%%%%%%%%%%%%%%%%%%%%%%%%%%%%%

\quad

We label the edges of the $(1,1)$-tangle presentation of the knot $8_9$
as in Figure~4.

\vspace{4mm}

%%%%%%%%%%%%%%%%%%%%%
\hspace*{43mm}
\includegraphics[scale=0.8]{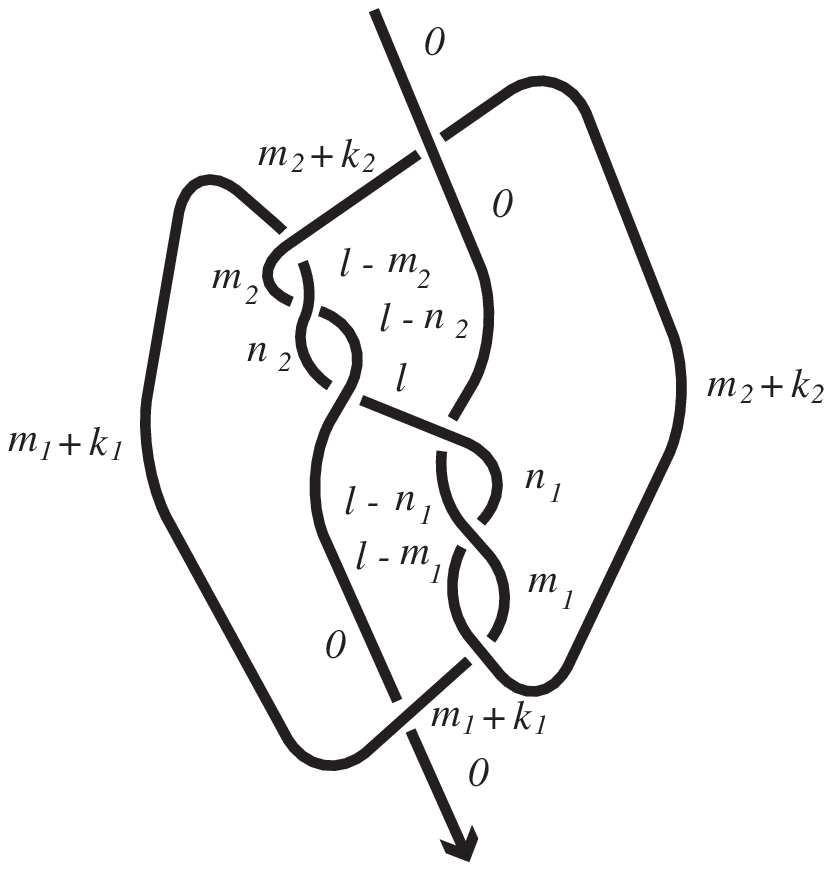}
%%%%%%%%%%%%%%%%%%%%%

\vspace{2mm}

\begin{center}
  Figure 4
\end{center}

\vspace{2mm}

We obtain the following formula of the colored Jones polynomial of the knot
$8_9$, where we put $l=m_1+m_2+k_1+k_2$ and use the formula
\eqref{eq:knot}:
\begin{align*}
  J_N (8_9) =
    & \sum_{\substack{0 \le l,m_1,m_2,n_1,n_2 \le N-1  \\
                      m_1+n_1,m_2+n_2 \le l  \\
                      m_1+m_2 \le l}}
      \left|
        \frac{(q)_{l-m_1}(q)_{l}(q)_{l-m_2}}
             {(q)_{m_1}(q)_{m_2}(q)_{n_1}(q)_{n_2}}
      \right|^{2}
      \frac{(\bar q)_{l-n_1} (q)_{l-n_2}}
           {(q)_{l-m_1-n_1} (\bar q)_{l-m_2-n_2}}  \\
    & \times
      q^{(m_2-m_1)(l-m_1-m_2)+(n_2-n_1)(l-n_1-n_2)+m_2-m_1+n_2-n_1},
\end{align*}
which can be asymptotically approximated by
$$
  \idotsint
  \exp \frac{\sqrt{-1}}{2\gamma} V_{8_9}(x,y,z,u,v)
  \, dx \, dy \, dz \, du \, dv ,
$$
where $x$, $y$, $z$, $u$, and $v$ correspond to $q^{-l}$, $q^{m_1}$, 
$q^{m_2}$, $q^{n_1}$, and $q^{n_2}$ respectively, and
\begin{align*}
  & V_{8_9}(x,y,z,u,v) \\
  & \hspace*{5mm}
    = - \li (x y) + \li \left( \frac{1}{x y} \right)
      - \li (x z) + \li \left(\frac 1 {x z}\right)
      - \li (x u) + \li \left( \frac{1}{x v} \right)  \\[2mm]
  & \hspace*{9mm}
      - \li (x) + \li \left( \frac{1}{x} \right)
      - \li (y) + \li \left( \frac{1}{y} \right)
      - \li (z) + \li \left( \frac{1}{z} \right)  \\[2mm]
  & \hspace*{9mm}
      - \li (u) + \li \left( \frac{1}{u} \right)
      - \li (v) + \li \left( \frac{1}{v} \right)
      + \li (x z v) -\li \left( \frac{1}{x y u} \right)  \\[2mm]
  & \hspace*{9mm}
    - \log \frac{y}{z} \log (xzv)
    - \log \frac{u}{v} \log (xyu). 
\end{align*}

Consequently we have
\begin{align*}
  & - \operatorname{Im}V_{8_9}(x_0,y_0,z_0,u_0,v_0) = 7.5881802 \ldots,
    \\[0.5mm]
  & \operatorname{Re} V_{8_9}(x_0,y_0,z_0,u_0,v_0) = 0
\end{align*}
for %the stationary point 
%$(x_0,y_0,z_0,u_0,v_0)=(0.7366011609-0.6763273835\sqrt{-1},
%                        0.4472176075-0.1647027124\sqrt{-1},
%                        1.968989044-0.7251455025\sqrt{-1},
%                        0.3859112582-0.0202712198\sqrt{-1},
%                        2.584139126-0.1357401508\sqrt{-1})$
\begin{align*}
x_0&=0.7366011609-0.6763273835\sqrt{-1},\\
y_0&=0.4472176075-0.1647027124\sqrt{-1},\\
z_0&=1.968989044-0.7251455025\sqrt{-1},\\
u_0&=0.3859112582-0.0202712198\sqrt{-1},\\
v_0&=2.584139126-0.1357401508\sqrt{-1}\\
\end{align*}
%\noindent                         
satisfying
\begin{gather*}
  \operatorname{Im}V_{8_9}(x_0,y_0,z_0,u_0,v_0)<0,  \quad
  \arg x_0 + \arg y_0+\arg u_0 \le 2\pi,  \\[0.5mm]
  \arg x_0 + \arg z_0+\arg v_0 \le 2\pi,  \quad
  \arg x_0 + \arg u_0+\arg v_0 \le 2\pi.
\end{gather*}
It follows from the calculation by SnapPea that
$$
  \exp \frac{\sqrt{-1}}{2\gamma} V_{8_9}(x_0,y_0,z_0,u_0,v_0)
  = \exp \frac{\vol(8_9)+\sqrt{-1} \CS(8_9)}{2\gamma},
$$
up to digits shown above.

%%%%%%%%%%%%%%%%%%%%%%%%%%%%%%%%%%%%%%%%%%%%%%%%%%%%%%%%%%%%%%%%%%%%%
\section{Knot $8_{20}$}
%%%%%%%%%%%%%%%%%%%%%%%%%%%%%%%%%%%%%%%%%%%%%%%%%%%%%%%%%%%%%%%%%%%%%

In this section, we  discuss a relation between the asymptotic behavior of 
the colored Jones polynomial and the Chern--Simons invariant for the  knot 
$8_{20}$.
We label each edge in the diagram of the knot in Figure~5. 

\vspace{2mm}

%%%%%%%%%%%%%%%%%%%%%
\begin{center}
  \includegraphics[scale=0.7]{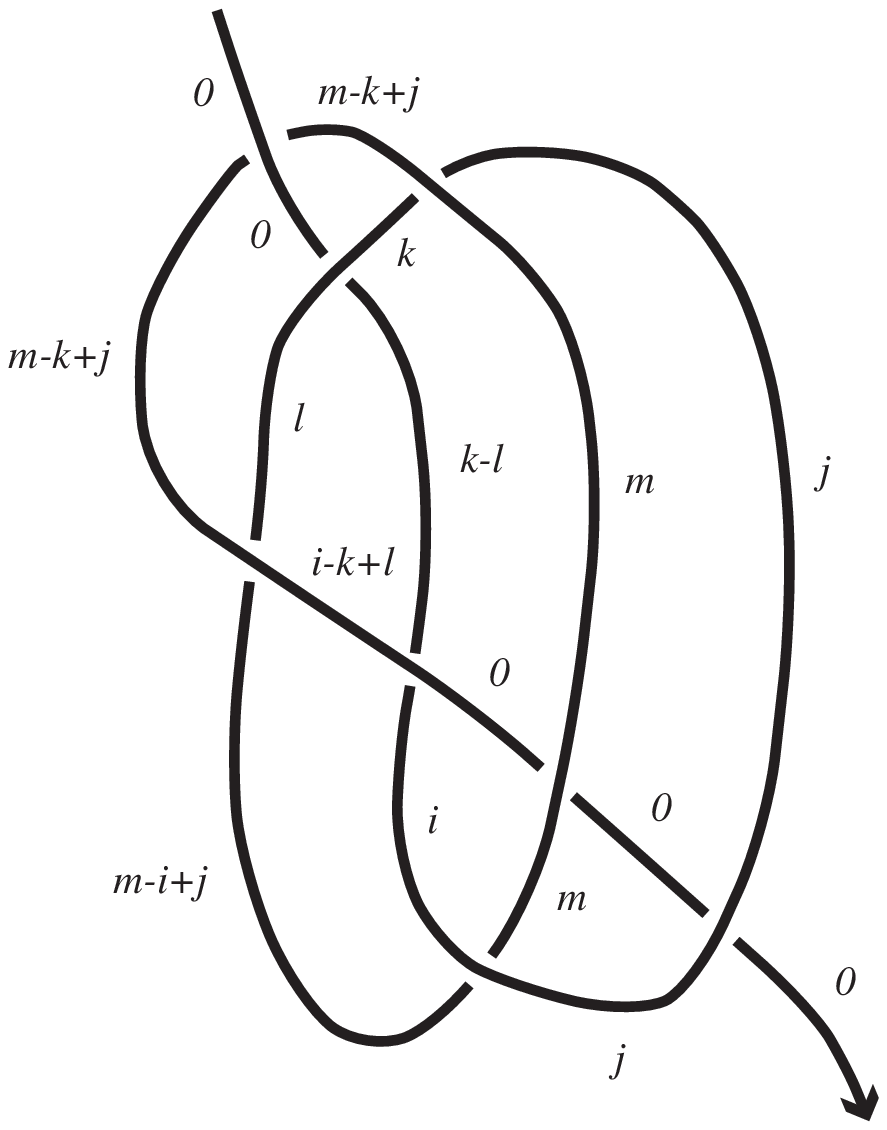}
\end{center}
%%%%%%%%%%%%%%%%%%%%%

\vspace{2mm}

\begin{center}
  Figure 5
\end{center}

\vspace{2mm}

The $N$-colored Jones polynomial of the knot $8_{20}$ is given by
\begin{equation}
  \sum_{\substack{j,l \leq k \leq i+l \leq j+m \\
                  j \leq i   \\
                  0 \le i,j,k,l,m \le N-1}}
 \frac{\{(\bar q)_{i} (q)_{k} (\bar q)_{m}\}^2}
      {\{ (\bar q)_{j} (q)_{l} \}^2
      (q)_{k-l}
       (\bar q)_{i-k+l} (\bar q)_{j+m-i-l} (q)_{i-j} (q)_{k-j}}
 q^{k+m+im+km-il}, \label{eq:knot820}
\end{equation}
which can be rewritten in the integral
\begin{equation*}
 \idotsint
 \exp \frac{\sqrt{-1}}{2\gamma} V_{8_{20}}(x,y,z,u,v)
 \,dx\,dy\,dz\,du\,dv
\end{equation*}
with 
\begin{align*}
  V_{8_{20}}(x,y,z,u,v)
    & = -2\li (x) + 2 \li \left(\frac 1 y \right) + 2 \li (z) 
        - 2\li \left(\frac 1 u\right)
        - 2\li \left(\frac 1 x\right)\\
    & \phantom{= \:}
      - \li \left(\frac 1 {xy}\right)
      - \li \left(\frac z {y}\right)  
      - \li (zu) +\li (xzu)  
      + \li \left(\frac 1 {xyuv}\right) \\
    & \phantom{= \:}
      + \log x \log u + \log x \log v -\log z \log v
      +\frac {\pi^2} 2.
\end{align*}
Here $x$, $y$, $z$, $u$, and $v$ correspond to $q^{-i}$, $q^j$, $q^{k}$, 
$q^{-l}$, and $q^{m}$ respectively. 

 Stationary points are solutions to partial differential equations
\begin{equation*}
  \frac{\partial V_{8_{20}}}{\partial x}
  = \frac{\partial V_{8_{20}}}{\partial y}
  = \frac{\partial V_{8_{20}}}{\partial z}
  = \frac{\partial V_{8_{20}}}{\partial u}
  = \frac{\partial V_{8_{20}}}{\partial v}
  = 0.
\end{equation*}
From these equations, we have the following system of algebraic equations
\begin{gather*}
  \left(1-x\right)^{2} \left(1-\frac{1}{xyuv}\right)uv
  =
  \left(1-\frac{1}{xy}\right) (1-xzu),  \\[2mm]
  \left(1-\frac 1 {xy}\right) \left(1-\frac{z}{y}\right) 
  =
  \left(1-\frac{1}{y}\right)^2 \left(1-\frac{1}{xyuv}\right),  \\[2mm]
  \left(1-z\right)^2 \left(1-xzu\right) v
  =
  \left(1-zu\right) \left(1-\frac{z}{y}\right) ,  \\[2mm]
  \left(1-zu\right)  \left(1-\frac{1}{xyuv}\right) x
  =
  \left(1-\frac{1}{u}\right)^2 \left(1-xzu\right),  \\[2mm]
  \left(1-\frac{1}{v}\right)^2 z
  =
  \left(1-\frac{1}{xyuv}\right) x.
\end{gather*}
Using MAPLE V, we get a stationary point 
$(x_0,y_0,z_0,u_0,v_0)$ 
%$(x_0,y_0,z_0,u_0,v_0)=(2.8785996+2.6574080\sqrt{-1},
%                          \infty,
%                          -0.4425377-0.45447889\sqrt{-1},
%                          0.3542198-0.0218067\sqrt{-1},
%                          0.1458832-0.3399257\sqrt{-1})$
%\noindent                           
which satisfies the conditions
\begin{equation*}
  \arg \frac{1}{u_0} \le \arg z_0,\quad
  \arg z_0 \le  \arg \frac{1}{x_0}+\arg \frac{1}{u_0}
\end{equation*}
from the range in the summation in \eqref{eq:knot820},
and 
\begin{equation*}
  \operatorname{Im}V_{8_{20}}(x_0,y_0,z_0,u_0,v_0)<0,
\end{equation*}
where $\operatorname{Im}$ denotes the imaginary part.
Note that the range of \eqref{eq:knot820} can be read as
\begin{gather*}
  \arg\frac{1}{u} \le \arg z \le \arg\frac{1}{x} +\arg\frac{1}{u}, \quad
  \arg\frac{1}{x} +\arg \frac {1} {y} +\arg \frac 1 {u}  \le  \arg{v},  \\[2mm]
  0  \le \arg\frac{1}{x}+ \arg \frac {1} {y} , \quad
  0\le \arg \frac{1}{x}, \arg z , \, \arg\frac{1}{u}, \, \arg v\le 2\pi.
\end{gather*}
To put it concretely,
\begin{align*}
x_0& =2.878599677+2.657408013\sqrt{-1},\\
y_0& =                          \infty,\\
z_0& =                          -0.4425377456-0.4544788919\sqrt{-1},\\
u_0& =                         0.3542198353-0.02180673815\sqrt{-1},\\
v_0& =                         0.1458832937-0.3399257634\sqrt{-1}.\\
\end{align*}

Then we obtain
\begin{align*}
  & - \operatorname{Im}V_{8_{20}}(x_0,y_0,z_0,u_0,v_0)
    = 4.1249032 \ldots ,  \\[1mm]
  & -\frac{\operatorname{Re} V_{8_{20}}(x_0,y_0,z_0,u_0,v_0) +\pi^2}
         {2\pi^2}
    = 0.1033634 \ldots .
\end{align*}
Applying values of $\vol(8_{20})$ and $\CS(8_{20})$ given by SnapPea 
\cite{Weeks:SnapPea}, we see that the following equation holds up to digits
shown above.
$$
  \exp \frac{\sqrt{-1}}{2\gamma} V_{8_{20}}(x_0,y_0,z_0,u_0,v_0)
  = \exp \frac {\vol(8_{20})+\sqrt{-1} \CS(8_{20})} {2\gamma}.
$$
Note that $\CS(8_{20})$ is defined modulo $\pi^2$. %$1/2$.

%%%%%%%%%%%%%%%%%%%%%%%%%%%%%%%%%%%%%%%%%%%%%%%%%%%%%%%%%%%%%%%%%%%%%
\section{Whitehead link}
%%%%%%%%%%%%%%%%%%%%%%%%%%%%%%%%%%%%%%%%%%%%%%%%%%%%%%%%%%%%%%%%%%%%%

\quad

For the final example, we calculate the limit of the colored Jones
polynomial of the Whitehead link given by \eqref{eq:Wh}, which can be
changed to the formula
\begin{equation*}
  J_N(L) = \sum_{\substack{0 \le i,j,k \le N-1 \\ k \le i,j}}
           \frac{\{(\bar q)_i (\bar q)_j\}^2} 
                {(q)_k^4 (\bar q)_{i-k} (\bar q)_{j-k}} q^{-(N-1)N/2}.
\end{equation*}

This can be asymptotically approximated by 
$$
  \int \! \! \! \int \! \! \! \int
  \exp \frac{\sqrt{-1}}{2\gamma} V_L(x,y,z) \, dx \, dy \, dz,
$$
where
$$
  V_L(x,y,z) = - 2 \li \left( \frac{1}{x} \right)
           - 2 \li \left( \frac{1}{y} \right)
           - 4 \li (z)
           + \li \left( \frac{z}{x} \right)
           +\li \left( \frac{z}{y} \right)
           + \pi^2,
$$
and $x$, $y$, and $z$ correspond to $q^i$, $q^j$, and $q^k$ respectively.
For a stationary point $(x_0, y_0, z_0) = (\infty, \infty, 1+\sqrt{-1})$, 
we obtain
\begin{align*}
  & - \operatorname{Im} V_L(x_0,y_0,z_0) = 3.663862 \ldots,  \\[1mm]
  & -\frac{\operatorname{Re} V_{L}(x_0,y_0,z_0)}{2\pi^2}
    = -0.1250000 \ldots.
\end{align*}
Since these values agree with SnapPea, the equation
$$
  \exp \frac{\sqrt{-1}}{2\gamma} V_{L}(x_0,y_0,z_0)
  = \exp \frac{\vol(L)+\sqrt{-1} \CS(L)}{2\gamma}
$$
holds up to digits shown above.

%%%%%%%%%%%%%%%%%%%%%%%%%%%%%%%%%%%%%%%%%%%%%%%%%%%%%%%%%%%%%%%%%%%%%
\section{Topological Chern--Simons invariant and some examples}
%%%%%%%%%%%%%%%%%%%%%%%%%%%%%%%%%%%%%%%%%%%%%%%%%%%%%%%%%%%%%%%%%%%%%

\quad

We propose a topological definition of the Chern--Simons invariant
for links.

 For a link $L$, if there exists the limit 
  \begin{equation*}
    2\pi\operatorname{Im} \lim_{N\to\infty} \log \frac {J_{N+1}(L)}{J_{N}(L)}
    \mod \pi^2, %\frac{1}{2}
  \end{equation*}
then we denote it by  $\CS_{\text{TOP}}(L)$ and 
call it the {\em topological Chern--Simons invariant} of $L$.

%%%%%
Let us give some numerical exapmles. 

For the knot $5_2$, we list some values of 
$(N, 2 \pi \log (J_{N+1}(5_2)/J_{N}(5_2)))$ by Pari-Gp in the following.
\begin{eqnarray*}
%{} & N & 2\pi \log (J_{N+1}(L)/J_{N}(L))\\
{}  & (40 , 3.058223721261842722613885956 - 3.022924613281720287391974968 \I) \\{} & (50 , 3.013081508530188353573854822 - 3.023340368517507069134855780 \I) \\
{} & (60 , 2.982744318753580696821772299 - 3.023574042878935429645720640 \I) \\
{} & (70 , 2.960955404961739170749114151 - 3.023717381786374852930574631 \I)\\
{} & (80 , 2.944548269170450112446966301 - 3.023811574968472287718611711 \I) \\
{} & (100 , 2.921483906108228993018469212 - 3.023923719027833555669502480 \I) \\{} & (120 , 2.906046421388666000282542398 - 3.023985374930307234443986632 \I) \\{} & (150 , 2.890559881907537128372001511 - 3.024036295143969179028770901 \I) \\{} & (200 , 2.875024234226941620327156350 - 3.024076266558545340852410631 \I)\\
{} & (250 , 2.865679250969538531562099056 - 3.024094905811349375139149331 \I) \\\end{eqnarray*}
By fitting the above data to  quadratic functions 
on $1/N$,
we can obtain the limit value 
$$2.82813-3.02414 \I$$
of $2\pi \log (J_{N+1}(5_2)/J_{N}(5_2))$ as 
$N \to \infty$  numerically, 
which  agrees with the value 
$$2.8281220 -3.02412837 \I $$ 
 by SnapPea.
 We display our data graphically in Figure 6 and Figure 7, wihich 
  help us to see the limit.

\quad

\hspace*{13mm}
%%%%%%%%%%%%%%%%%%%%%
\includegraphics[scale=1]{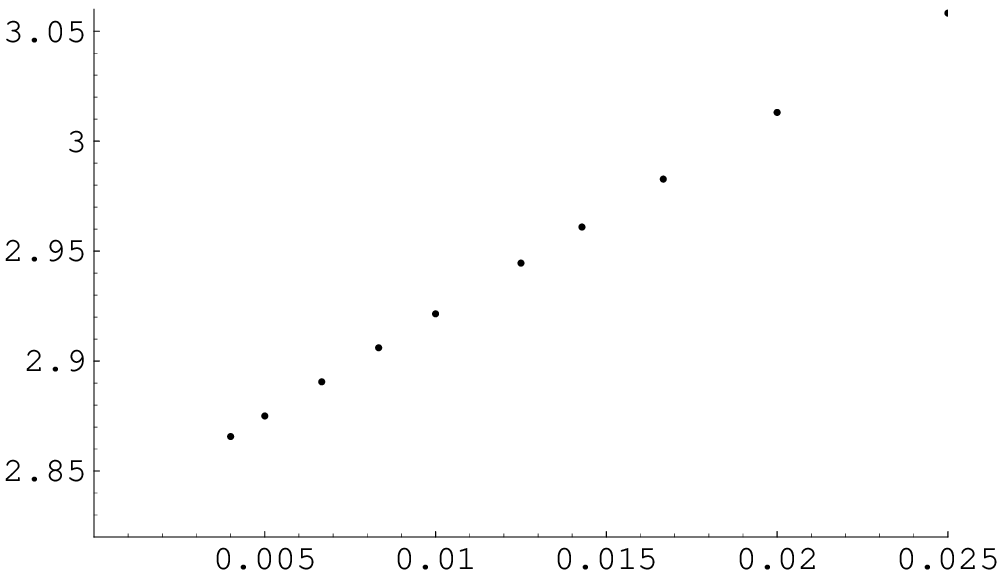}
%%%%%%%%%%%%%%%%%%%%%

\begin{center}
Figure 6. Dots indicate 
$(1/N , 2\pi \operatorname{Re}\log (J_{N+1}(5_2)/J_{N}(5_2)))$ 
for $N=40,50,60,70,80,100,120,150,200,250$. The origin corresponds to 
$(0, 2.82)$. 
\end{center}

\quad

\hspace*{13mm}
%%%%%%%%%%%%%%%%%%%%%
\includegraphics[scale=1]{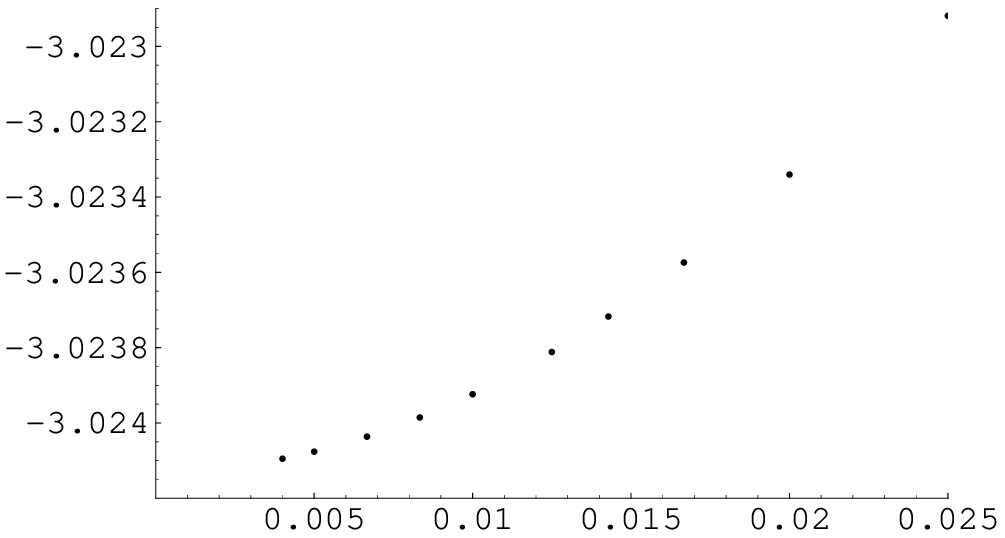}
%%%%%%%%%%%%%%%%%%%%%

\begin{center}
Figure 7. Dots indicate 
$(1/N , 2\pi \operatorname{Im}\log (J_{N+1}(5_2)/J_{N}(5_2)))$ 
for $N=40,50,60,70,80,100,120,150,200,250$. The origin corresponds to 
$(0, -3.0242)$.
\end{center}

\quad

Similarly, for the Whitehead link $L$, we illustrate our numerical check
 in Table 1, 
Figure 8, and Figure 9.

\quad

Table 1. $(N , 2\pi \log (J_{N+1}(L)/J_{N}(L)))$ for the Whitehead link $L$
\begin{eqnarray*}
{} & (40 , 3.892920359101811097809525583 + 2.457483997330866045812504703 \I) \\
{} & (50 , 3.848161466402914225154530180 + 2.461039474018016569869745301 \I) \\
{} & (60 , 3.818029013349499312708236153 + 2.462976748675980254703390855 \I) \\
{} & (70 , 3.796362501209537691078944556 + 2.464147191795881614582476451 \I) \\
{} & (80 , 3.780034327560022195082015385 + 2.464907923404764622274395868 \I) \\
{} & (100 , 3.757062258985477857247991239 + 2.465803785962819679236327339 \I) \\{} & (120 , 3.741674608179023673159144258 + 2.466291085896660260688606142 \I) \\{} & (150 , 3.726228649726558590507828429 + 2.466690204011030007962113880 \I) \\\end{eqnarray*}

\hspace*{13mm}
%%%%%%%%%%%%%%%%%%%%%
\includegraphics[scale=1]{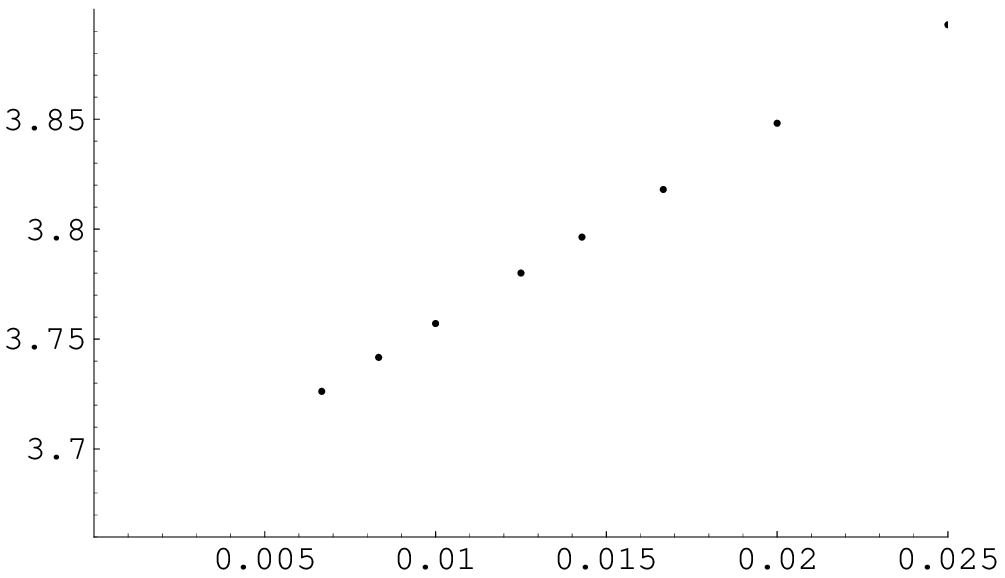}
%%%%%%%%%%%%%%%%%%%%%

\begin{center}
Figure 8. Dots indicate 
$(1/N , 2\pi \operatorname{Re}\log (J_{N+1}(L)/J_{N}(L)))$ 
 for  $N=40,50,60,70,80,100,120,150$. 
The origin corresponds to $(0, 3.66)$. 
\end{center}

\hspace*{13mm}
%%%%%%%%%%%%%%%%%%%%%
\includegraphics[scale=1]{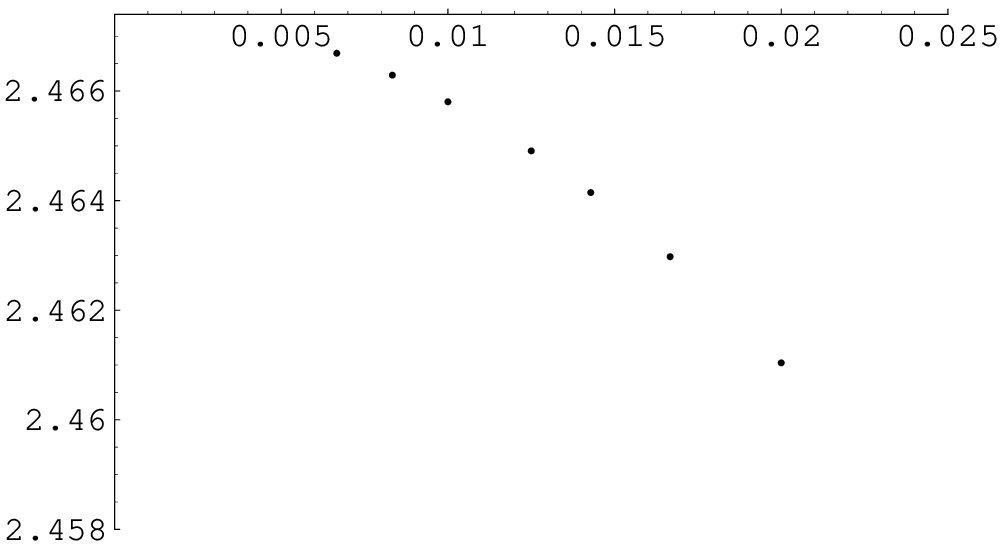}
%%%%%%%%%%%%%%%%%%%%%

\begin{center}
Figure 9. Dots indicate 
$(1/N , 2\pi \operatorname{Im}\log (J_{N+1}(L)/J_{N}(L)))$ 
for $N=40,50,60,70,80,100,120,150$. 
The origin corresponds to $(0, 2.4674)$. 
\end{center}

\quad

Fitting, we get the numerical limit value
 $3.66386+2.46742 \I$
of $2\pi \log (J_{N+1}(L)/J_{N}(L))$ as 
$N \to \infty$, 
which  agrees with our result in the section 6.

%%%%%

%%%%%%%%%%%%%%%%%%%%%%%%%%%%%%%%%%%%%%%%%%%%%%%%%%%%%%%%%%%%%%%%%%%%%
\section{Conclusion}
%%%%%%%%%%%%%%%%%%%%%%%%%%%%%%%%%%%%%%%%%%%%%%%%%%%%%%%%%%%%%%%%%%%%%

\quad

We have shown the following by concrete calculations.
\begin{obs}\label{obs:volume+CS}
  Let $L$ be one of the hyperbolic knots $6_3$, $8_9$, and $8_{20}$,
  or the Whitehead link. 
  Following Kashaev's way, we approximate the colored Jones polynomial
  $J_N(L)$ of $L$  asymptotically by 
  $$
    \int \cdots \int
    \exp \frac{N\sqrt {-1}}{2\pi} V_L({\bf x}) d{\bf x}.
  $$
  Then there exists a stationary point ${\bf x}_0$ of $V_L$ such that
  the formula
  \begin{equation*}
    \exp \frac{N\sqrt{-1}}{2\pi} V_{L}({\bf x}_0)
    = \exp \frac{N}{2\pi} (\vol(L)+\sqrt{-1} \CS(L))
  \end{equation*}
  holds up to $6$ digits.
\end{obs}

\begin{conj}[Complexification of Kashaev's conjecture]
Let $L$ be a hyperbolic link. 
Then, it  holds that 
\begin{equation*}
  \vol(L)
  =
  2 \pi \lim_{N\to \infty} \frac{\log\left|\langle L \rangle_N\right|}{N}
\end{equation*}
with $\vol(L)$ the hyperbolic volume of the complement of $L$. 
Moreover, there exists the topological Chern--Simons invariant 
$\CS_{\text{TOP}}(L)$  of $L$
\begin{equation*}
   \CS_{\text{TOP}}(L)= 
   2\pi\operatorname{Im} \lim_{N\to\infty} \log \frac {J_{N+1}(L)}{J_{N}(L)}
    \mod \pi^2, %\frac{1}{2},
\end{equation*} 
and $\CS_{\text{TOP}}(L)$ equals to $CS(L)$ modulo $\pi^2$.
Here $\CS(L)$ is the Chern--Simons invariant of $L$ 
\cite{Chern/Simons:ANNMA274,Meyerhoff:density}.
Note that the complement of $L$ is a hyperbolic manifold with cusps.
\end{conj}

We note that Observation~\ref{obs:volume+CS} also holds for the knots
$4_1$, $5_2$ and $6_1$ by calculating Kashaev's examples in
\cite{Kashaev:LETMP97} using MAPLE V and SnapPea.

Therefore we conclude that the complexified Kashaev's conjecture is true,
up to several digits, up to choices of contours when we change summations
into integrals, and up to choices of saddle (stationary) points when we
approximate integrals by the saddle point method, for the six hyperbolic
knots above and for the Whitehead link.

Note that if the complexified Kashaev's conjecture is true then the
topological Chern--Simons invariant of a hyperbolic link coincides with its
Chern--Simons invariant associated with the hyperbolic metric.
Moreover if the volume conjecture is true then the colored Jones polynomial
would give both the simplicial volume and the topological Chern--Simons
invariant for any knot.

%%%%%%%%%%%%%%%%%%%%%%%%%%%%%%%%%%%%%%%%%%%%%%%%%%%%%%%%%%%%%%%%%%%%%
\bibliography{mrabbrev,mmoty-bib}
\bibliographystyle{amsplain}
%%%%%%%%%%%%%%%%%%%%%%%%%%%%%%%%%%%%%%%%%%%%%%%%%%%%%%%%%%%%%%%%%%%%%

\end{document}